\documentstyle[12pt]{article}
\renewcommand{\baselinestretch}{1.3}

\makeatletter
\newcommand{\single}{\let\CS=\@currsize\renewcommand{\baselinestretch}{1}\tiny\CS}
\newcommand{\singles}{\let\CS=\@currsize\renewcommand{\baselinestretch}{1.3}\tiny\CS}
\newcommand{\oneanda}{\let\CS=\@currsize\renewcommand{\baselinestretch}{1.2}\tiny\CS}
\newcommand{\doubles}{\let\CS=\@currsize\renewcommand{\baselinestretch}{1.5}\tiny\CS}
\newcommand{\tree}{\let\CS=\@currsize\renewcommand{\baselinestretch}{1.5}\tiny\CS}
\newcommand{\four}{\let\CS=\@currsize\renewcommand{\baselinestretch}{2}\tiny\CS}

\newcommand{\ncom}{\newcommand}
\ncom{\bq}{\begin{equation}}
\ncom{\eq}{\end{equation}}
\ncom{\beqn}{\begin{eqnarray*}}
\ncom{\eeqn}{\end{eqnarray*}}
\ncom{\beq}{\begin{eqnarray}}
\ncom{\eeq}{\end{eqnarray}}
\ncom{\been}{\begin{enumerate}}
\ncom{\eeen}{\end{enumerate}}
\ncom{\nno}{\nonumber}
\ncom{\hs}{\mbox{\hspace{.25cm}}}
\ncom{\rar}{\rightarrow}
\ncom{\lrar}{\longrightarrow}
\ncom{\Rar}{\Rightarrow}
\ncom{\noin}{\noindent}
\newtheorem{thm}{Theorem}[section]
\newtheorem{lemma}[thm]{Lemma}
\newtheorem{cor}[thm]{Corollary}
\newtheorem{pro}[thm]{Proposition}
\newtheorem{example}[thm]{Example}
\newtheorem{remark}[thm]{Remark}

\ncom{\bt}{\begin{thm}}
\ncom{\et}{\end{thm}}
\ncom{\bl}{\begin{lemma}}
\ncom{\el}{\end{lemma}}
\ncom{\bco}{\begin{cor}}
\ncom{\eco}{\end{cor}}
\ncom{\bp}{\begin{pro}}
\ncom{\ep}{\end{pro}}
\ncom{\bex}{\begin{example}}
\ncom{\eex}{\end{example}}
\ncom{\brm}{\begin{remark}}
\ncom{\erm}{\end{remark}}
\ncom{\comx}{I\!\!\!\!C}
\ncom{\zee}{$Z\!\!\!\!Z$}
\ncom{\ze}{Z\!\!\!\!Z}
\ncom{\Q}{$I\!\!\!\!Q$}
\ncom{\p}{I\!\!P}
\ncom{\al}{\alpha}
\ncom{\be}{\beta}
\ncom{\f}{\frac}
\ncom{\ga}{\gamma}
\ncom{\bib}{\bibitem}
\ncom{\pf}{{\bf Proof: }}
\ncom{\sta}{\stackrel}
\ncom{\cA}{{\cal A}}
\ncom{\cG}{{\cal G}}
\ncom{\cI}{{\cal I}}
\ncom{\cO}{{\cal O}}
\ncom{\cV}{{\cal V}}
\ncom{\cW}{{\cal W}}
\ncom{\cK}{{\cal K}}
\ncom{\cM}{{\cal M}}
\ncom{\cL}{{\cal L}}
\ncom{\cZ}{{\cal Z}}
\title{ Abelian Varieties into $2g+1$-dim. Linear Systems}
\author{Jaya N.Iyer}

\begin{document} 
\maketitle
\footnotetext{Mathematics Classification Number: 14C20, 14B05, 14E25.}
\begin{abstract} We show that polarisations of type $(1,...,1,2g+2)$ on
  $g$-dimensional abelian varieties are $\it{never}$ very ample, if
  $g\geq 3$. This disproves a conjecture of Debarre, Hulek and
  Spandaw. We also give a criterion for non-embeddings of abelian
  varieties into $2g+1$-dimensional linear systems.
\end{abstract}

\section{Introduction}
Let $L$ be an ample line bundle of type $\delta$ on an abelian variety
$A$, of dimension $g$. Classical results of Lefschetz ($n\geq 3$) and Ohbuchi ($n=2$)
imply very ampleness of $L^n$, if $|L|$ has no fixed divisor when $n=2$.
Suppose $L$ is an ample line bundle of type $(1,...,1,d)$ on $A$.
When $g=2$, Ramanan ( see [4]) has shown that if $d\geq 5$ and the
abelian surface does not contain elliptic curves then $L$ is very
ample. When $g\geq 3$, Debarre, Hulek and Spandaw ( see [3],
Corollary 25, p. 201) have shown the following.

\bt
Let $(A,L)$ be a generic polarized abelian variety of dimension $g$
and type $(1,...,1,d)$. For $d>2^g$, the line bundle $L$ is very ample.
\et
They further conjecture that if $d\geq 2g+2$, then the line bundle $L$
is very ample, ( see [3], Conjecture 4, p. 184).
In particular, when $g=3$ and $d\geq 8$, their results ( for $d\geq 9$)
and conjecture ( for $d=8$) imply that $L$ is very ample.

The results due to Barth ( [1])and Van de Ven ( [5])show
\bt
For $g\geq 3$, no abelian variety $A_g$ can be embedded in $\p^d$,
  for $d\leq 2g$.
\et
In particular, it implies that line bundles of type $(1,...,1,d)$,
$d\leq 2g+1$, are never very ample.

We show
\bt
Suppose $L$ is an ample line buundle of type $(1,...,1,d)$ on an
abelian variety $A$, of dimension $g$. If $g\geq 3$ and $d\leq 2g+2$,
then $L$ is never very ample.
\et
 
This disproves the conjecture of Debarre et.al when $d=2g+2$ and
gives a different proof of 1.2, for morphisms into the complete linear
system $|L|$. The proof of 1.3 also indicates the type of singularities of
the image in $|L|$.

Now any abelian variety $A$ of dimension $g$ can be embedded in a
projective space of dimension $2g+1$. 

Consider a morphism $A\lrar |V|$, where $dim|V|=2g+1$.
Suppose the involution $i:A\lrar
A$, $a\mapsto -a$ lifts to an involution on the vector space $V$,
hence on the linear system $|V|$, ( such a situation will arise, for
example, if $A$ is embedded by a symmetric line bundle into its 
complete linear system, of dimension greater than $2g+1$.
 One may then project the abelian
variety from a vertex which is invariant for the involution $i$, to a
projective space of dimension $2g+1$ and the involution $i$ will then descend down
to this projection ). 

Then we show
\bt
Suppose there is a morphism $A\sta{\phi}{\lrar} |V|$, with $dim|V|=2g+1$ and the
involution $i$ acting on the vector space $V$. If $degree
\phi(A)>2^{2g}$ and $dim V_+\neq dim V_-$, then the morphism $\phi$ is never an embedding, for
all $g\geq 1$. Here $V_+$ and $V_-$ denote the $\pm 1$-eigenspaces of $V$,
for the involution $i$.
\et

\section{Proof of 1.3.}
Consider a pair $(A,L)$, as in 1.3.

We may assume $L$ is an ample line bundle of characteristic $0$ on
$A$. Then $L$ is symmetric , i.e. there is an isomorphism
$L\simeq i^*L$, for the involution $i:A\lrar A,\,a\mapsto
-a$. This induces an involution on the vector space $H^0(L)$, also
denoted as $i$. Let $H^0(L)^+$ and $H^0(L)^-$ denote the $+1$ and
$-1$-eigenspaces of $H^0(L)$, for the involution $i$ and $h^0(L)^+$ and
$h^0(L)^-$ denote their respective dimensions.
Choose a normalized isomorphism $\psi:L\simeq i^*L$, i.e. the fibre
map $\psi(0):L(0)\lrar L(0)$ is $+1$.

Let $A_2$ denote the set of torsion $2$ points of $A$. If $a\in A_2$
then $\psi(a):L(a)\lrar L(a)$ is either $+1$ or $-1$.

Let $$A_2^+=\{a\in A_2: \psi(a)= +1\}$$
and $$A_2^-=\{a\in A_2:\psi(a)=-1\}$$
and $Card(A_2^+)$ and $Card(A_2^-)$ denote their respective
cardinalities.

Consider the associated morphism $A\sta{\phi_L}{\lrar }\p H^0(L)$ and
let $$\p _+= \p \{s=0 : s\in H^0(L)^-\}$$ and $$\p _-=\p \{ s=0 : s\in
H^0(L)^+\}.$$
Then the involution $i$ acts trivially on the subspaces $\p _+$ and
$\p _-$, of $\p H^0(L)$. Moreover, $\phi_L(A_2^+)\subset \p _+$ and
$\phi_L(A_2^-)\subset \p _-$.

\bl
If $a\in A_2^+$, then the intersection of the image $\phi_L(A)$ and
$\p _+$ is transversal at the point $\phi_L(a)$.
\el
\pf The action of the involution $i$ at the tangent space, $T_{A,a}$, at
$a$, is $-1$. If the intersection of $\phi_L(A)$ with $\p _+$ is not
transversal at $\phi_L(a)$, then $\phi_{L*}(T_{A,a})$ intersects
$\p_+$, giving a $i$-fixed non-trivial subspace of $T_{A,a}$, which is
not true.  
( This argument was given by M.Gross.)
 $\Box$

Let $Z=\phi_L(A)\cap \p _+$, in $\p H^0(L)$. Then
$\phi_L(A_2^+)\subset Z$.
Suppose $dimZ>0$. Since the involution $i$ acts trivially on $Z$, the
morphism $\phi_L$ restricts on $\phi_L^{-1}(Z)\lrar Z$, as a morphism
of degree at least $2$, with its Galois group containing $<i>$.
If $dimZ=0$, then by 2.1, the points of $\phi_L(A_2^+)$ have
multiplicity $1$ in $Z$. Let $r=deg Z-Card(A_2^+)$. Then there are
$\frac{r}{2} $-points on $\phi_L(A)$ on which the involution $i$ acts
trivially, i.e. there are $\frac{r}{2}$-pairs $(a,-a),\,a\in A-A_2$,
which are identified transversally by $\phi_L$. By $K(L)$-invariance
of the image $\phi_L(A)$, there are more such pairs.
\brm
If $dimZ>0$ or $r>0$, then $L$ is not very ample.
\erm
$\bf{Case\,1:}$
$d=2m $ and $m\leq g+1$.

By [2], 4.6.6, $h^0(L)^+=m+1$ and $h^0(L)^-=m-1$.

Hence $dim\p _+=m$ and $dim \p_-=m-2$.

a) If $m<g+1$, then $dimZ\geq g+m-2m+1>0$.

b) If $m=g+1$. By Riemann-Roch, $deg\phi_L(A)=(2g+2).g!$.
If $dimZ=0$ then since $\p_+$ and $\phi_L(A)$ have complementary
dimensions in $\p H^0(L)$, $degZ=(2g+2).g!$.

Now by [2], Exercise 4.12 b)-Remark 4.7.7,
$$Card(A_2^+)\leq 2^{2g-(g-1)-1}(2^{g-1}+1)$$
$$=2^g(2^{g-1}+1).$$
Since $g\geq 3$, $r\geq(2g+2).g!-2^g(2^{g-1}+1)>0$.

Hence by 2.2, $L$ is not very ample.

$\bf{Case\,2:}$
$d=2m-1$ and $m\leq g+1$.

Then $h^0(L)^+=m$ and $h^0(L)^-=m-1$. Hence $dim\p _+=m-1$ and
$dim\p_-=m-2$.

a) If $m<g+1$, then $dimZ\geq g+m+1-2m>0$.

b) If $m=g+1$, as in $\bf{Case\,1}$, $deg\phi_L(A)=(2g+1)g!$, and $\p
_+$ and $\phi_L(A)$ have complementary dimension in $\p H^0(L)$. Hence
if $dimZ=0$, then $degZ=(2g+1)g!$. Also, in this case, $Card(A_2^+)\leq
2^{g-1}(2^g+1)$.

Since $g\geq 3$, $r\geq(2g+1)g!-2^{g-1}(2^g+1)>0$. Hence by 2.2, $L$ is
not very ample.   $\Box$

\section{Morphisms into $i$-invariant linear systems}

$\bf{Proof\,of\,1.4}$: 
Consider the morphism $A\sta{\phi}{\lrar}|V|$, with the involution $i$
acting on the vector space $V$.
Let
 $$\p_+= \p \{s=0: s\in V_-\}$$ and 
 $$\p_-= \p \{s=0: s\in V_+\},$$
where $V_+$ and $V_-$ denote the $+1$ and $-1$-eigenspaces of the 
vector space $V$, for the involution $i$. Let $d=degree \phi(A)$.

Now $dim\p_+>g$ or $dim\p_+< g$ or $dim\p_+=g$.

Case 1: $dim\p_+> g$

Consider the intersection $Z=\p_+\cap \phi(A)$.

Then $dimZ \geq g+ g+1 -2g-1\geq 0$.

As in Proof of 1.3, 
if $dimZ>0$, then the restricted morphism $\phi^{-1}(Z)\lrar Z$ is of
degree at least $2$, since $i$ acts trivially on $Z$. 
Suppose $dim Z=0$. Then the
intersection of $\phi(A)$ and $\p_+$ is transversal at the image of
torsion $2$ points of $A$, by 2.1. Since
$Card(A_2)=2^{2g}$ and $degree(\phi(A))>2^{2g}$, there are pairs $\{a,-a\}$
on $A$ which get identified transversally by the morphism $\phi$.

Case 2: $dim\p_+< g$.

In this situation, $dim \p_-> g$ and we can repeat the above argument.
  
Hence $\phi$ is never an embedding.
 $\Box$

\brm When $dim V_+= dim V_-$, the morphism $\phi$ need not identify some pair
of points $\{a,-a\}$, in the linear system $|V|$. For example, consider a symmetric line
bundle $L$, of type $(1,1,9)$, on a generic abelian threefold $A$. Then $L$ is very ample and
$dim H^0(L)_+=5$ and $dim H^0(L)_-=4$. Hence $dim \p_+=4$ and $dim \p_-=3$. Consider the scroll $S_A=\cup_{a\in A} l_{a,-a
}$, where $l_{a,-a}$ is the line joining the points $a$ and $-a$, in $ |L|$.
Then the line $l_{a,-a}$ is invariant for the involution $i$ and has two fixed points,
one of them in $\p_+$ and the other in $\p_-$. Hence $S_A$ intersects $\p_+$ in at most
a $3-$dimensional subset. Now we can project from a point of $\p_+$, outside this subset, and
the projection will have the fixed spaces of $i$ to be equidimensional. Also, by the choice of
the point of projection, there are no pairs $\{a,-a\}$ identified in the projection.

\erm

$\it{Acknowledgement}$: We thank M.Gross for giving the argument in
 2.1 and O.Debarre for pointing out a mistake in an earlier version dealing with the
 case when $dim V_+= dim V_-$, in 1.4.
 We also thank French Ministry of National Education, Research and Technology, for their
 support.

\begin{thebibliography}{99}

\bib [1]{1}  Barth, W.:{\em Transplanting cohomology classes in
  complex-projective space}, Amer.J. of Math. $\bf{92}$, 951-967, (1970).

\bib [2]{2}  Birkenhake, Ch., Lange, H. : {\em Complex abelian varieties},   
 Springer-Verlag, Berlin, (1992).

\bib [3]{3}  Debarre, O., Hulek, K., Spandaw, J. : { \em Very ample
linear systems on abelian varieties}, Math. Ann. $\bf{ 300}$, 181-202, (1994).

\bib [4]{4}  Ramanan, S.: {\em Ample Divisors on Abelian Surfaces},
Proc. London Math. Soc. (3), $\bf{51}$, 231-245, (1985).

\bib [5]{5} Van de Ven  : {\em On the embeddings of abelian
  varieties in projective spaces}, Ann.Mat.Pura Appl.(4),
  $\bf{103}$, 127-129, (1975).

\end {thebibliography}

\noindent{Institut de Mathematiques,\\ Case 247, Univ.Paris-6, \\4, Place
  Jussieu, \\75252, Paris Cedex 05, France.\\ Email: iyer@math.jussieu.fr }

\end{document}